\documentclass[journal]{IEEEtran}
\ifCLASSINFOpdf
\else
   \usepackage[dvips]{graphicx}
\fi
\usepackage{url}

\begin{document}
%
\title{Mean field variational framework for integer optimization}
%
%
%

\author{Arturo~Berrones,
        Jon\'as~Velasco,
        and~Juan~Banda,
\thanks{A. Berrones, J. Velasco and J. Banda are with 
Universidad Aut\'onoma de Nuevo Le\'on,
Facultad de Ingenier\' \i a Mec\'anica y El\'ectrica
Posgrado en Ingenier\' \i a de Sistemas,
AP 126, Cd. Universitaria, 
San Nicol\'as de los Garza, NL 66450, M\'exico
e-mail: arturo@yalma.fime.uanl.mx.}}

\maketitle

\begin{abstract}
A principled method to obtain approximate solutions of general constrained integer optimization 
problems is introduced. The approach is based on the calculation of a mean field probability distribution for
the decision variables which is consistent with the objective function and the constraints. The original discrete task
is in this way transformed into a continuous variational problem. 
In the context of particular problem classes 
at small and medium sizes, the mean field results
are comparable to those of standard specialized methods, while at large sized instances
is capable to find feasible solutions in computation times for which
standard approaches can't find any valuable result. 
The mean field variational framework remains valid for
widely diverse problem structures so it represents a promising paradigm for large dimensional
nonlinear combinatorial optimization tasks.
\end{abstract}


%
\IEEEpeerreviewmaketitle

\section{Introduction}
%
%
%
%
\IEEEPARstart{B}{esides} of its practical importance and of the progress made in the solution of integer linear programs in the last decades 
\cite{hemmecke},
integer optimization is still a very challenging subject. Practical methods for 
nonlinear problems are rare and mostly limited to small to medium sized problems \cite{ambrosio,hemmecke,buchheim}, 
even in special cases such
as convex integer programming \cite{buchheim,bonami,aydin}. 
Heuristic approaches consequently play a major role in the field, 
either hybridized with deterministic algorithms \cite{belotti,dark}, or
like stochastic approximate solution strategies on their own right \cite{hartmann,hard,anneal,ga,eva,yung,yang,xiao,narayan}. 

An essential ingredient in any heuristic is a rule by which candidate solutions are generated. In this contribution we introduce
a very general framework to construct probabilistic solution generators that are consistent with the underlying integer optimization task.

 

\section{Mean Field Framework}

The method is principled and firmly based on the fundamental rules of probability. 
More precisely, any possible solution for an integer optimization problem can be viewed like a random string 
drawn from some probability distribution. It is desirable that such a distribution generates deviates that are tipically in the feasibility region
and close to the optima. To proceed, 
consider the following class of optimization problems,
\begin{eqnarray} \label{opt}
\min f(\vec{x}) \quad s.t. \\ \nonumber
g_{k}(\vec{x}) \leq 0, h_{l} (\vec{x}) = 0,
\end{eqnarray}
where $\vec{x}$ is a vector of binary decision variables, $g_{k (k=1,..., K)}$ are inequality constraint functions
and $h_{l (l=1,..., L)}$ are equality constraints. The optimization task (\ref{opt}) can in principle be represented by a potential
function $V(\vec{x})$ which includes the objective and the constraints. A probability distribution can be associated to such a potential by
the transformation \cite{anneal,metropolis},
\begin{eqnarray} \label{boltzmann}
P(\vec{x}) = \frac{1}{Z} \exp \left (-\frac{V}{kT} \right ),
\end{eqnarray}
where $Z$ is a normalization factor (or partition function) and $kT$ is a constant. The Eq. (\ref{boltzmann}) gives the maximum
entropy distribution which is consistent with the condition $\left < V \right >_{P} = \int V P d{\vec{x}}$ \cite{jaynes}.
However, $P$ is in general intractable. Moreover, in our setup is not even known, because the explicit definition of $V$ would require
the knowledge of suitable ``barrier'' terms that exactly represent the constraints.
Is therefore proposed the following mean field probabilistic model for the decision variables,
\begin{eqnarray} \label{pmodel1}
Q(\vec{x}) = \prod_{i=1}^{N} p(x_i),
\end{eqnarray}
Mean field techniques, which have first emerged in statistical mechanics \cite{parisi}, have been 
already successfully applied to discover fundamental features of combinatorial problems and valuable
solution strategies, although focused on particular combinatorial problem classes \cite{martin}. 
Our purpose in this contribution is to develop a practical mean field framework to find good candidate
solutions to linear and nonlinear integer problems in the constrained situation (\ref{opt}).
The most general form for the independent marginals is,
\begin{eqnarray} \label{pmodel2}
p(x_i) = 1 + (2m_{i} - 1)x_{i} - m_{i}.
\end{eqnarray}
The $m$'s are continuous mean field parameters, $m \in [0,1]$. These parameters can be selected by
the minimization of the Kullback-Leibler divergence between distributions $Q$ and $P$ \cite{opper},
\begin{eqnarray}
D_{KL} (Q || P) =  \left < \ln Q \right > - \left < \ln P \right >,
\end{eqnarray}
where the brackets represent averages with respect to the tractable distribution $Q$. Introducing the entropy $S_Q = -kT \left < \ln Q \right >$,
is obtained the variational problem $\min F_Q$, where
\begin{eqnarray}
F_Q = \frac{1}{kT} \left [\left < V \right > - S_Q \right]
\end{eqnarray}
is the variational ``free energy'' of the distribution $Q$ \cite{opper}.
Without loss of generality, the constant $kT$ can be set $kT = 1$. In first instance we consider the class of combinatorial optimization
problems in which all the involved functions (objective and constraints) are polynomial, 
e. g. $f(\vec{x}) = a_o + \sum_i b_i x_i + \sum_i \sum_j c_{i,j} x_i x_j + \sum_i \sum_j \sum_r q_{i,j,r} x_i x_j x_r  + ... $. 
In such case 
$\left < f(\vec{x}) \right > = f(\left < \vec{x} \right >)$, $\left < g(\vec{x}) \right > = g(\left < \vec{x} \right >)$
and $\left < h(\vec{x}) \right > = h(\left < \vec{x} \right >)$. Therefore, the continuous relaxation of the problem (\ref{opt})
is equivalent to its average under the mean field distribution,
\begin{eqnarray} \label{mopt}
\min f(\vec{m}) \quad s.t. \\ \nonumber
g_{k}(\vec{m}) \leq 0, h_{l} (\vec{m}) = 0.
\end{eqnarray}
An expression for $\left < V \right >$ can be constructed in terms of the Lagrangian,
\begin{eqnarray} \label{pot}
\left < V \right >  = \mathcal{L} = f(\vec{m}) + \sum_{l} \lambda_l h_{l} (\vec{m}) + \sum_{k} \mu_k g_{k}(\vec{m}),
\end{eqnarray}
where the constants $\lambda_l$ and $\mu_k \geq 0$ are the Karush-Kuhn-Tucker (KKT) multipliers \cite{chong}. 
The entropy of $Q$, on the other hand, is given by,
\begin{eqnarray}
S_Q = - \sum_{i} \left [ (1-m_i)\ln(1-m_i) + m_i \ln m_i \right],
\end{eqnarray}
so the variational problem for the $m$'s is written like,
\begin{eqnarray} \label{pmodel3}
\min F_Q (\vec{m}) =
\min  \{
f(\vec{m}) + \sum_{l} \lambda_l h_{l} (\vec{m}) + \sum_{k} \mu_k g_{k}(\vec{m}) \\ \nonumber
+ \sum_{i} \left [ (1-m_i)\ln(1-m_i) + m_i \ln m_i \right]  \}.
\end{eqnarray}
Equations (\ref{pmodel1}), (\ref{pmodel2}) and (\ref{pmodel3}) give a general probabilistic model for combinatorial optimization problems with binary decision variables. 
Any continuous and differentiable nonlinearities in the objective or the constraints can be expanded in a Taylor series under to the condition $m_i < 1$ $\forall$ $i$.
Due to independence under the mean field, $\left < V \right >$ is therefore given by Eq. (\ref{pot}) for any problem (\ref{opt}), provided 
that the stated conditions are met.
Stationarity applied to $F_Q$ with respect to the mean field parameters reduce the variational 
problem to a set of self-consistency equations for $\vec{m}$,
\begin{eqnarray}\label{gsol}
m_i = \frac{1}{1 + \exp[\partial_i \mathcal{L}_{\vec{\lambda},\vec{\mu}}(\vec{m})]}
\end{eqnarray}
The generation of valid solutions of (\ref{opt}) from
the mean field model 
can then be tackled by the following numerical scheme, 
\begin{enumerate}
\item Give initial values for the KKT multipliers.
\item Solve Eq. (\ref{gsol}).
\item Evaluate the objective and constraints, either rounding the mean field (MF)
solution or drawing a point from the resulting MF distribution.
Update the KKT multipliers 
and repeat from stage $2$ until suitable stoping criteria are met.
\end{enumerate}

\section{Results}

The probabilistic setup (\ref{pmodel1}), (\ref{pmodel2}), (\ref{pmodel3}) and (\ref{gsol})
is now tested on specific examples. We first consider here 
the classical Knapsack Problem (KP), 
\begin{eqnarray} \label{kp}
\min \quad - \vec{q} \cdot \vec{x} \quad s.t. \\ \nonumber
\vec{w} \cdot \vec{x} - d \leq 0,
\end{eqnarray}
where $d$ is the capacity of the knapsack, $\vec{q}$ are the gains and $\vec{w}$ the weights of a collection of $i = 1, ..., N$ objects.
The set of self-consistency equations for the mean field parameters is in this case
independent, with solution,
\begin{eqnarray}\label{solkp}
m_i = \frac{1}{1 + \exp(- q_i + \mu w_i)},
\end{eqnarray}
where $\mu$ is the KKT multiplier associated to the single constraint. When this constraint is inactive, $\mu = 0$ and $m_i \in \left [\frac{1}{2}, 1 \right )$.
In the case in which the constraint is active, $\mu \geq 0$.
The substitution of Eq. (\ref{solkp}) into the complementary slackness condition $\vec{w} \cdot \vec{m} - d = 0$ gives,
\begin{eqnarray}\label{slackness}
&&\sum_{i=1}^{N}(w_i - d)\exp(q_i - \mu w_i) \\ \nonumber
&+& \sum_{i=1}^{N} (w_i + w_{i+1}- d) 
\exp \left[ (q_i + q_{i+1}) 
- \mu (w_i + w_{i+1}) \right ] \\ \nonumber
&+& \sum_{i=1}^{N} (w_i + w_{i+1} + w_{i+2} - d) \\ \nonumber
&&\exp \left[ (q_i + q_{i+1} + q_{i+2}) - \mu (w_i + w_{i+1} + w_{i+2}) \right ] \\ \nonumber
&+& ... + \sum_{i=1}^{N} (\tau \overline{w} + \epsilon_{\tau} - d) 
\exp \left[\tau \overline{q} + \epsilon^{'}_{\tau} - \mu (\tau \overline{w} + \epsilon_{\tau}) \right ] \\ \nonumber
&+& ... +
\left[ \left( \sum_{i=1}^{N}w_i \right ) - d \right ] \exp \left( \sum_{i=1}^{N} q_i - \mu \sum_{i=1}^{N} w_i \right ) 
= d,
\end{eqnarray}
where $q_{i+j} = w_{i+j} = 0$ if $(i+j) > N$.
The overlines represent averages over the distributions from which the instance parameters (weights and gains) are drawn.
The index $\tau = 1, ..., N$ counts the terms in the left side of Eq. (\ref{slackness}). The finite sums over 
the instance parameters that appear in the expression are rewritten in terms of the estimators,
\begin{eqnarray}
\frac{1}{\tau} \sum_{t=1}^{\tau} w_{i+t} = \overline{w} + \epsilon_{\tau}, \\ \nonumber
\frac{1}{\tau} \sum_{t=1}^{\tau} q_{i+t} = \overline{q} + \epsilon^{'}_{\tau},
\end{eqnarray}
where the $\epsilon$'s are estimation errors, which in general go to zero as $\tau \to \infty$. For the problem to be
feasible, $\overline{w} << d$ if $N >> 1$. This imply that for large $N$,
the terms of the left side of Eq. (\ref{slackness}) with
small $\tau$ are typically negative, and should be cancelled out. Therefore,
\begin{eqnarray}\label{slackness2}
\alpha (N \overline{w} - d)\exp[N (\overline{q} - \mu \overline{w})] \approx d, \quad \alpha \geq 1,
\end{eqnarray}
from which it follows that,
\begin{eqnarray}\label{mu}
\mu \approx \frac{\overline{q}}{\overline{w}} - \frac{1}{\alpha N} \ln \left ( \frac{d}{N\overline{w} - d} \right ).
\end{eqnarray}
In the limit $N \to \infty$, $\mu = \frac{\overline{q}}{\overline{w}}$. This result gives an important insight: for active constraints, the multipliers  
are such that the distribution (\ref{pmodel2}) actually represents the ``competition'' between the 
objective and the constraints. Good starting points of the search scheme for the multipliers can be obtained on this basis. For KP, we have
implemented the following procedure,
\begin{enumerate}
\item The multiplier is initially set like $\mu = \frac{\overline{q}}{\overline{w}}$.
\item Evaluate Eq. (\ref{solkp}).
\item Update $\mu$ by the minimization of $(\vec{w} \cdot \vec{m} - d)^2$.
Evaluate $\vec{w} \cdot \vec{m} - d = \epsilon$. Given a predefined tolerance $tol$,
if $|\epsilon| > tol$, then go to step $2$ using a new initial value for the multiplier drawn
at random from a neighboor around $\frac{\overline{q}}{\overline{w}}$. 
Else, end.
\end{enumerate}
At each iteration the resulting mean field parameters are rounded, the feasibility is checked
and the objective is evaluated, keeping track of the best feasible solution. The random 
initial values for the multiplier are taken from an interval of $10$\% around the $\mu$ value
associated to the current best rounded solution at each iteration.
Some numerical results and comparisons are presented in Table 1. 
In the reported experiments, $tol = 0.0001$. A maximum of $1000$ iterations is allowed.
The scheme is terminated if $|\epsilon| \leq tol$ or the total number
of permitted iterations is completed.
The best found rounded feasible solution and its objective value is reported.
The instances were created with the 
Pisinger generator \footnote{\url{http://www.diku.dk/~pisinger/generator.c}} \cite{kp},
under the conditions of strong linear correlations between $q_i$ and $w_i$, being $w_i$ randomly distributed
in the interval $[1, 1000]$.
It has been argued that for these type of instances 
the integer solution
is usually far from the continuous relaxation solution \cite{kp}.
The exact solutions have been found using the branch and bound algorithm provided 
by Cplex \footnote{\url{http://www-01.ibm.com/software/commerce/optimization/cplex-optimizer/index.html}}, 
which is a
widely accepted state of the art standard for linear and quadratic integer problems. 
A general purpose Genetic Algorithm (GA) is also tested on the instances.
The Cplex version is Ibm Ilog Cplex 12.4. For the GA the R package Genalg 
\footnote{\url{http://cran.r-project.org/web/packages/genalg/}} 
version 0.1.1 was used,
with all the defaults except a rule which stops GA evolution if there is no change in best objective value after
$150$ generations or a time limit of $5000$ seconds of CPU time is exceeded.
The experiments were run on an Intel I7-3820 processor with 4 cores, 3.60 GHz and 32 GB RAM.
For each problem dimension, $10$ instances are drawn from the Pisinger generator, 
running each instance $10$ times from different initial conditions. 
Averages and 
standard deviations of the objective values and computation times for the MF and GA algorithms are reported.
For the considered instances, Cplex found the exact solutions in CPU times $< 1$ second.
Altough GA and MF are comparable at small problem sizes, for large problems
the mean field heuristic
finds higher quality solutions at a computation time that is orders of magnitude below of GA. 

\noindent The effectivity of the formulation in nonlinear cases has been tested on the Quadratic Knapsack Problem (QKP),
which is stated as follows \cite{pisinger2},
\begin{eqnarray}
\min \quad - \vec{x}^{t}\mathcal{Q}\vec{x} \quad s.t. \\ \nonumber
\vec{w} \cdot \vec{x} - d \leq 0,
\end{eqnarray}
where $\mathcal{Q}$ is a symmetric matrix with coefficients $q_{i,j} \geq 0$ $\forall$ $i,j$. 
QKP has a graph-theoretic interpretation in terms of the {\em Clique} problem \cite{pisinger2}.
Moreover, QKP is in general NP-hard in the strong sense, meaning that 
admits no polynomial-time approximation scheme unless P = NP \cite{pisinger2}.
The mean field Lagrangian reads,
\begin{eqnarray}
\mathcal{L} =  - \vec{m}^{t}\mathcal{Q}\vec{m} + \mu (\vec{w} \cdot \vec{m} - d),
\end{eqnarray}
from which
\begin{eqnarray}
\partial_i \mathcal{L}_{\mu} =  -2 q_{i,i} m_i - \sum_{j \neq i} q_{i,j} m_j + \mu w_i.
\end{eqnarray}
The initial value of the multiplier is taken like,
\begin{eqnarray}
\mu = \frac{1}{N \bar{w}}\left [ 2 \sum_{i} q_{i,i} + \sum_i \sum_{j\neq i} q_{i,j} \right ].
\end{eqnarray}
Due to nonlinearity, the algebraic system (\ref{gsol}) is now coupled.
At fixed $\mu$, it can be efficiently solved by iterating the equations (\ref{gsol}).
Besides this, the numerical procedure is almost the same followed for the KP example.
More precisely, the step $2$ of the procedure now involves the iteration of the mean field
equations (\ref{gsol}) from a starting random initial $\vec{m}$. In our experiments, we have used only
one iterating step.
Results in instances of $200$ variables
taken from the benchmark provided by Billionet and Soutif 
\footnote{\url{http://cedric.cnam.fr/~soutif/QKP/QKP.html}}
\cite{bill} 
are reported in Table 2.
Even in these relatively small instances Cplex is unable to
find optimality certificates in short times. 
Therefore we now include the branch and bound method in the comparison, 
focusing on solution quality and computation times. The equipment is the same as that
used for the KP experiments.
A maximum CPU time of $100$ seconds is allowed for the mean field and the branch and bound
algorithms, while for the GA is imposed a limit of $5000$ seconds. 
The averages and standard deviations of the ratio between the best value found with respect the
known optimum is reported. Each instance has been run $10$ times, while the benchmark includes $10$
examples of each of the quadratic coefficients matrix densities of $25$\%, $50$\% and $100$\%. 
Therefore the statistics is computed from
a total of $300$ runs per method.
Both the branch and bound and the mean field methods have an average
solution quality above $99$ \% with respect to the known exact optimum, while the GA shows
considerable inferior performance. Larger instances are studied using the 
QKP generator provided by Pisinger 
\footnote{\url{http://www.diku.dk/~pisinger/testqkp.c}}.
Comparisons are in this case only between the mean field and branch and bound by computing the ratio
between the best solutions found by MF and Cplex. 
A total of $10$ instances is generated for each problem size, 
with matrix density of $100$\%. Each instance was ran $10$ times from different
initial conditions, giving a total statistics of $100$ runs for each problem size per method.
For problem sizes
of $1000$ the performance between both methods appear to be similar, however at size $2000$ 
Cplex is unable to find integer solutions in some cases. These are discarded until the $10$ examples are completed.
For sizes of $5000$ and above,
the branch and bound is unable to find any integer solution in the $100$ seconds computation time limit. 
In fact, by further tests it have turned out that Cplex is unable to give integer solutions for these cases
even after several days of CPU time in our equipment.
We have
tested the mean field up to instances of size $20000$, always giving integer feasible solutions, although with
an increasing variance of the best objectives values. 
Figure 1 shows that the variance is however reduced and the 
overall solution quality incremented when larger 
computation times are allowed, which indicates that the mean field procedure scales robustly for 
large problem dimensions.

\section{Conclusion}

Under mild assumptions on the statistical properties of the linear KP, an explicit mean field solution has been found
in the infinite size limit, expressed by the equations (\ref{solkp}) and (\ref{mu}). 
In the experiments presented in
Table 1, the ratio between the initial value of $\mu$ given by $\mu = \frac{\overline{q}}{\overline{w}}$ 
and the final value
obtained after the iterations of the numerical scheme has been found to be $1 \pm 0.003$. 
This kind of ``thermodynamic limit'' might exist in other linear constrained problems as well.
The nonlinear constrained case studied here also points out
in a similar direction. In all the considered problems, the heuristically proposed initial multiplier
always give an integer valid solution, despite that this is hard to achieve by a branch and bound
search from random initial conditions. 

The mean field framework presented in this contribution appears to be highly competitive with respect to standard
approaches. 
For the test problems considered, there exist other more specialized deterministic and 
stochastic methods besides the GA and branch and bound used here for comparisons. However, the available
literature on these other algorithms is 
mainly focused on small to medium problem sizes \cite{yung,yang,xiao,narayan}.
We have used therefore two popular stochastic and deterministic methods as a way to asses the scalability of our
framework with respect to what is often done in practice.
Moreover, our interest is the development of a general paradigm. The implementation for
almost arbitrary nonlinear problem structures should follow essentially the same steps used here for the 
QKP example. Clearly, the application of our framework 
to different integer optimization problem classes is one of the research lines opened by the present work. 
Other is the design of novel heuristics based on the 
mean field or its hybridization with existing algorithms. 
Also could be valuable to explore corrections to the mean field equations,
by the use of cavity, replica or related approaches \cite{parisi,opper}.

\begin{table*}[h] 
\caption{Classical knapsack problem. 
Statistics of the ratio of the best values obtained
by GA and MF with respect to the known 
optimum (found by the branch and bound algorithm implemented in Cplex) 
and of the CPU times per run (in seconds) is computed from
$10$ instances for each problem size, running each instance $10$ times from random and 
independent
initial conditions for both methods. 
On these and 
the following reported experiments, all the runs were performed on the same equipment under
the same controlled conditions.}
\begin{tabular*}{\hsize}{@{\extracolsep{\fill}}ccccc}
\hline
\multicolumn1c{$N$}&\multicolumn1c{[(GA best) / (optimum)] (\%)}&\multicolumn1c{GA time}&  
\multicolumn1c{[(MF best) / (optimum)](\%)}&\multicolumn1c{MF time}
\cr
\hline
$100$ & $95.44 \pm 0.66$ & $24.88 \pm 5.37$ & $92.35 \pm 5.68$ & $0.11 \pm 0.05$ \cr
$1000$ & $91.34 \pm 3.35$ & $224.01 \pm 65.24$ & $94.32 \pm 2.23$ & $4.192 \pm 3.66$ \cr
$10000$ & $90.43 \pm 3.43$ & $1098.71 \pm 170.05$ & $96.51 \pm 0.62$ & $48.19 \pm 27.73$ \cr 
\hline
\end{tabular*}
\end{table*}

\begin{table*}[h] 
\caption{Quadratic knapsack problem. The problem size with $N = 200$ decision variables from the
Billionet and Soutif benchmark is considered. 
A total of $100$ seconds of CPU time is allowed for MF and Cplex, while a $5000$ seconds
of CPU time limit is imposed to the GA. 
The benchmark consists of $10$ instances
of each of the quadratic coefficients matrix densities of $25$\%, $50$\% and $100$\%.
Statistics were collected by performing $10$ independent runs per method on each instance. 
The ratio of the best values with respect to the
known optima is therefore reported from a total of $300$ runs per method .}
\begin{tabular*}{\hsize}{@{\extracolsep{\fill}}ccccc}
\hline
\multicolumn1c{[(GA best) / (optimum)] (\%)}&\multicolumn1c{[(MF best) / (optimum)] (\%)}&
\multicolumn1c{[(Cplex best) / (optimum)] (\%)}
\cr
\hline
$87.4_{-6.4}^{+6.5}$ & $99.4_{-0.2}^{+0.5}$ & $99.9_{-0.1}^{+0.0}$\cr 
\hline
\end{tabular*}
\end{table*}

\begin{table*}[h] 
\caption{Quadratic knapsack problem. Large problem dimensions are considered by drawing them from the
Pisinger generator. No exact optima are known in these cases. The ratio between the best found
values of MF and Cplex are reported from $10$ different instances ran $10$ times by both methods, imposing a 
limit of $100$ seconds of CPU time. At size $N = 2000$, Cplex is unable to find integer solutions for some cases.
These have been discarded until $10$ instances solvable by Cplex are completed. At size $N = 5000$, Cplex gave none
integer solution.}
\begin{tabular*}{\hsize}{@{\extracolsep{\fill}}ccccc}
\hline
\multicolumn1c{$N$}&\multicolumn1c{[(MF best) / (Cplex best)] (\%)}
\cr
\hline
$500$ & $99.5_{-2.0}^{+0.4}$ \cr
$1000$ & $100.0_{-0.0}^{+0.4}$ \cr
$2000$ & $101.2_{-1.6}^{+4.2}$  \cr 
$5000$ & $\infty$ \cr
\hline
\end{tabular*}
\end{table*}

\vskip 0.5cm
\begin{figure*}[h]
\begin{center}
{\includegraphics[width=1.0\textwidth]{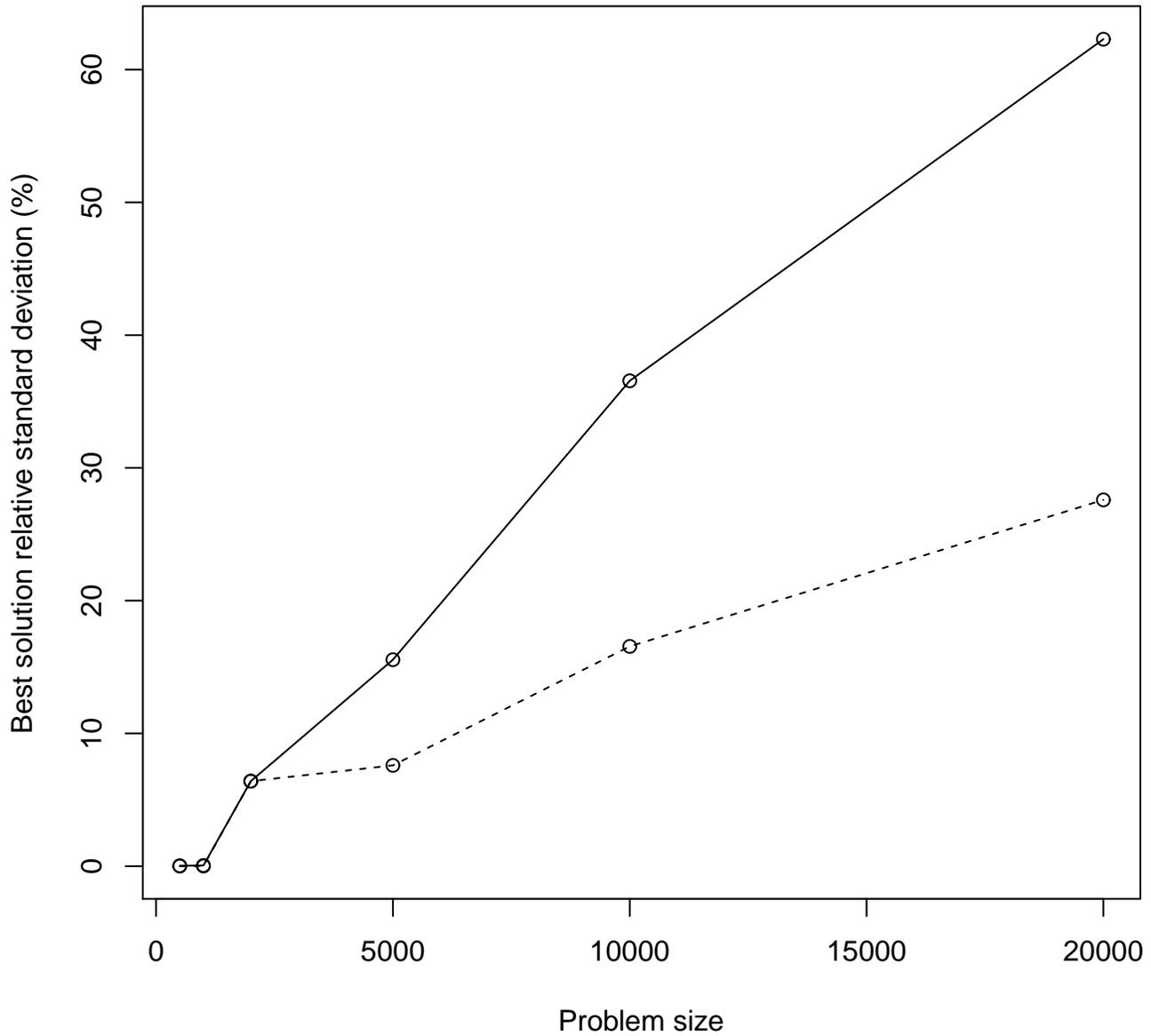}}
\caption{The scalability of the MF procedure is tested by analyzing the relative standard deviation
(the ratio of the standard deviation of the best values with respect to the average of the best values, in 
percent units) of the MF results at different sizes of the QKP. 
The solid curve reports the behavior with
$100$ seconds of CPU time while the dashed curve gives the relative standard deviations with $200$ seconds of
CPU time. The curves show that integer feasible solutions can be found by MF in modest
computation times for problem sizes  
that are intractable by standard methods. The variability of MF decreses by
increasing computation times. The ratio between average best values at $200$ and $100$ seconds is
$1.13_{-0.1}^{+0.3}$. Therefore the figure is indicative of an overall increase in solution quality as a consequence of the
CPU time increment.
}\label{fig1}
\end{center}
\end{figure*}

\section*{Acknowledgment}

This work was partially supported by the National Council of Science and Technology of
Mexico under grant CONACYT CB-167651 and by the Autonomous University of Nuevo Leon 
support to research program under grant UANL-PAICYT IT795-11.

\ifCLASSOPTIONcaptionsoff
  \newpage
\fi

\end{document}